\documentclass[12pt]{amsart}
\usepackage{amsmath}

\def\CC{{\rm\kern.24em\vrule width.02em height1.4ex depth-.05ex\kern-.26em C}}
\def\PP{{\rm\kern.24em\vrule width.02em height1.4ex depth-.05ex\kern-.26em P}}
\def\ZZ{{\rm\kern.26em\vrule width.02em height0.5ex depth0ex\kern.04em\vrule width.02em height1.47ex depth-1ex\kern-.34em Z}}
\def\BB{{\rm\kern.24em\vrule width.02em height1.4ex depth-.05ex\kern-.26em B}}

\newcommand {\BOX} {\rule{2mm}{2mm} \bigskip}
\newcommand{\ra}{\rightarrow}

\newcommand {\ol} {\overline}

\newtheorem{lem}{Lemma}[section]
\newtheorem{theo}[lem]{Theorem}

\newtheorem{remark}[lem]{Remark}

\numberwithin{equation}{section}

\begin{document}
\title{Compositional roots of H{\'e}non maps}
\author{Gregery T. Buzzard and John Erik Forn\ae ss}
\date{}


\begin{abstract}
Let $H$ denote a composition of complex H\'enon maps in $\CC^2$.
In this paper we show that the only possible compositional roots of
$H$ are also compositions of H\'enon maps, and that $H$ can have
compositional roots of only finitely many distinct orders.
\end{abstract}

\renewcommand{\thefootnote}{}
\footnote{Research at MSRI is 
   supported in part by NSF grant DMS-9022140}

\maketitle
\let\cal\mathcal

\section{Introduction}

Following \cite{bs1}, we say that a generalized H{\'e}non map is a map
of the form
$$ H(z,w) = (w, p(w) - az), $$
where $p$ is a monic polynomial of degree $d \geq 2$ and $a \in
\CC-\{0\}$, and we let ${\cal G}$ denote the space of finite compositions
of such maps.  From \cite{fm}, we know that any polynomial
diffeomorphism of $\CC^2$ is conjugate either to one of the maps in
${\cal G}$ or to an elementary map which preserves each line of the
form $w = \text{const}$.  In \cite{bf}, we classified, up to
conjugacy, all polynomial diffeomorphisms which arise as the time-1 map
of a holomorphic vector field.  In particular, each of these maps is
an elementary map and has
compositional roots of all orders.  Moreover, in some cases, these
roots can be nonpolynomial.  See \cite{forst} for information about
such cases.  

\bigskip
In this paper we treat the question of the existence of compositional
roots for the remaining cases.  In particular, we show that any root
of a map in ${\cal G}$ must be a polynomial map and that any map in
${\cal G}$ can have roots of
only a finite number of distinct orders.  For the remaining elementary
maps which are not the time-1 map of a flow, we show that such maps
have roots of arbitrarily high order and nonpolynomial roots,
but that any root of such a map is conjugate to a polynomial
elementary map.

\section{Dynamical behavior and Green's functions}

Fix $H \in {\cal G}$.  Let $K^+ (= K^+(H))$ and
$K^- (= K^-(H))$ denote the
set of points $p$ in $\CC^2$ such that the orbit of $p$ under $H$ is
bounded with respect to forward or backward iteration, respectively.
Also, for $R>0$, define sets
\begin{align*}
V^- & := \{(z,w): |w|>R \text{ and } |w|>|z|\}, \\
V^+ & := \{(z,w): |z|>R \text{ and } |w|<|z|\},\\
V & := \{(z,w): |z| \leq R \text{ and } |w| \leq R\}.
\end{align*}
A simple argument shows that $K^\pm \subseteq V^\pm \cup V$ for $R$
sufficiently large.
Finally, we let $d$ denote the degree of $H$ as a polynomial map, and
define functions 
$$ G^\pm(z,w)  :=  \lim_{n \ra \infty} \frac{1}{d^n} \log^+\|H^{\pm
n}(z,w)\|. $$

\bigskip
It is immediate that $G^+$ is continuous and plurisubharmonic (psh)
on $\CC^2$, identically $0$ on $K^+$, and strictly positive and
pluriharmonic on $\CC^2 - K^+$.  Analogous statements are true with
$G^-$ and $K^-$ in place of $G^+$ and $K^+$.   See, for example,
\cite{bs1} for a more systematic discussion.

Note also that $G^\pm \circ H = d^{\pm 1} G^\pm$.  Moreover,
corollary~2.6 in \cite{bs1} implies the following.

\begin{lem}   \label{lemma:bounds}
There exist $R>0$, $C>0$ such that if $(z,w) \in \ol{V^- \cup V}$,
then 
$$ \log^+|w| - C \leq G^+(z,w) \leq \log^+|w| + C,$$
and if $(z,w) \in \ol{V^+ \cup V}$, then
$$ \log^+|z| - C \leq G^-(z,w) \leq \log^+|z| + C. $$
\end{lem}

Since $G^+(z_0,w)$ is subharmonic in $w$ for each fixed $z_0$, we see
that it is bounded from above in $\{|w|<|z_0|\}$ by its maximum on the
boundary of this disk, which is contained in $\ol{V^- \cup V}$.
Applying a similar argument to $G^-$ gives the following.

\begin{lem}   \label{lemma:max}
There exists $C>0$ such that for $(z,w) \in \CC^2$, 
$$ G^\pm (z,w) \leq \max\{\log^+|z|, \log^+|w|\} + C. $$
\end{lem}

\section{Compositional roots and Green's functions}

Fix $H \in {\cal G}$ and let $d$ be the degree of $H$.  In this
section we show that if $F^n = H$ in the sense of composition, then
$G^\pm \circ F = d^{\pm 1/n} G^\pm$.  First a simple lemma.

\begin{lem}
Suppose $F$ is an automorphism of $\CC^2$ and $F^n = H$.  Then
$K^+$ and $K^-$ are the same for $F$ as for $H$, and $F$ is a
diffeomorphism of $K^+$ and of $K^-$.
\end{lem}

{\bf Proof:}  Take $p \in K^+(H)$ and let ${\cal O}^+_H(p)$ denote the
forward orbit of $p$ under $H$.  Then $\cup_{j=0}^{n-1} F^j(\ol{{\cal
O}^+_H(p)})$ is compact, and ${\cal O}^+_F(p)$ is contained in
this set, hence is bounded.  Thus $K^+(H) \subseteq K^+(F)$.

If $p \not \in K^+(H)$, then the forward orbit is not bounded for $H =
F^n$, hence is not bounded for $F$.   Thus $K^+(H) = K^+(F)$, and a
similar argument applies to $K^-$. 

The fact that $F$ is a diffeomorphism of $K^\pm = K^\pm(F)$ is clear
from the definition of these sets.  $\BOX$  \bigskip

\begin{lem}   \label{lemma:mult}
Let $F$ be as in the previous lemma.  Then $G^\pm \circ F = d^{\pm 1/n} G^\pm$.
\end{lem}

{\bf Proof:}  Since $F$ is holomorphic on $\CC^2$ and preserves $K^+$,
we see that $G^+ \circ F$ is $0$ on $K^+$, plurisubharmonic and 
continuous on $\CC^2$, and strictly positive and pluriharmonic on
$\CC^2 - K^+$.  

Fix $z_0$ and define $g_{z_0}(w) := G^+ \circ F(z_0,w)$.  Note that $K^+
\cap (\{z_0\} \times \CC)$ is a compact set and that $g$ is harmonic
on the complement of this set.  Hence, outside a large disk, 
$g_{z_0}$ has a harmonic conjugate in a neighborhood of any point.  Using
analytic continuation in the exterior of this disk, we obtain a
harmonic conjugate with periods.  Hence for some $r>0$, some constant
$c_{z_0}$, and a real harmonic function $h_{z_0}$, we get a function 
$$ \phi_{z_0}(w) = g_{z_0}(w) - c_{z_0} \log|w| + i h_{z_0}(w)$$
which is holomorphic for $|w| > r$.

Since $g_{z_0} \geq 0$, we have $|\exp(-\phi_{z_0}(w))| \leq |w|^{c_{z_0}}$.
Hence $\exp(-\phi_{z_0}(w))$ has at most a pole at infinity, so we can write 
$$ \exp(-\phi_{z_0}(w)) = w^N \exp(f(w))$$
for some integer $N$ and some $f$ holomorphic in $|w|>r$ with a
removable singularity at infinity. 

Taking absolute value and $\log$, we get $g_{z_0}(w) - c_{z_0} \log|w|  = -N
\log|w| - \text{Re}(f(w))$.  Hence $g_{z_0}(w) = b_{z_0} \log|w| +O(1)$ in
$\{|w|>2r\}$, for some $b_{z_0}$.  Since $g_{z_0} \geq 0$, we have
$g_{z_0}(w) = b_{z_0} \log^+|w| + O(1)$ in $\CC$.  

\bigskip
We claim that $b_z$ is independent of $z$.  Note that $2 \pi
b_{z_0}$ is the period for the harmonic conjugate of $g_{z_0}$ in
$|w|>r$, and that $g_z(w)$ is pluriharmonic in $(z,w)$ near
$(z_0,w_0)$ for any $|w_0|>r$. 

Fix $|w_0|>r$.  We can construct the harmonic conjugate for $g$ in the
bidisk $\Delta(z_0;r_0) \times \Delta(w_0;r_0)$ for some $r_0$ small.  For
$w \in \Delta(w_0;r_0)$, we can use analytic continuation as above to
extend $g_z$ around a circle in $\{z\} \times \{|w|>r\}$.  Doing
this for each such $w$ gives a new harmonic conjugate for $g$ in the
bidisk, which must differ from the original by a constant.  Thus
$b_z = b_{z_0}$ for $z$ near $z_0$.  

\bigskip
Hence $g_{z_0}(w) = b \log^+|w| + O(1)$, and from lemma~\ref{lemma:bounds},
we see $G^+(z_0,w) =  \log^+|w| + O(1)$.  Thus $g_{z_0}(w) - b
G^+(z_0,w)$ is continuous for $w \in \CC$ and harmonic for $w$ such
that $(z_0, w) \notin K^+$, has a removable singularity at $\infty$,
and is $0$ for $w$ such that $(z_0, w) \in K^+$, which is a nonempty
set.  Hence $g_{z_0}(w) \equiv b G^+(z_0,w)$ for all $w \in \CC$.

Similarly, $G^+ \circ F(z,w) - b G^+(z,w) \equiv 0$ 
for all $(z,w) \in \Delta(z_0;r_0) \times \CC$.
Since this difference is pluriharmonic in
$\CC^2-K^+$, which is connected, it must be $0$ throughout $\CC^2 -
K^+$, hence throughout $\CC^2$ since both terms are $0$ on $K^+$.

Finally, $G^+ \circ F^n(z,w) = G^+ \circ H(z,w) = d
G^+(z,w)$, while induction with the above result shows that $G^+ \circ
F^n(z,w) = b^n G^+(z,w)$  Hence $b^n = d$, and $b>0$ since
$G^+\geq 0$.  This gives the lemma for $G^+$, and the same proof
applies to $G^-$.    $\BOX$ \bigskip

\section{Polynomial roots}

In the proof of the following theorem, we use the terminology and
results of \cite{fm}.  In particular, we use the fact that the group
of polynomial automorphisms of $\CC^2$ is the amalgamated product of
the group ${\cal A}$ of affine linear automorphisms and the group
${\cal E}$ of elementary automorphisms which preserve the set of lines
of the form $w = const$.  A {\em reduced word} is an automorphism of the
form $g_1 \cdots g_k$, $k \geq 1$, where each $g_k$ is in ${\cal A}$ or ${\cal E}$
but not in the intersection of these two groups and no two adjacent
$g_j$'s are in the same group.  We say that $k$ is the {\em length} of this
reduced word.  Also, we need to know that the identity cannot be
written as a reduced word.  

\begin{theo}
Suppose $H \in {\cal G}$ is a composition of generalized H{\'e}non
maps and $F$ is an automorphism of $\CC^2$ with $F^n = H$.  Then $F
\in {\cal G}$.
\end{theo}

{\bf Proof:}  Let $(z,w) \in \CC^2$, and let $F=(F_1,F_2)$.  If
$F(z,w) \in \ol{V^- \cup V}$, then from 
lemmas~\ref{lemma:bounds}, \ref{lemma:mult}, and \ref{lemma:max}, we
see
\begin{align*}
\log^+|F_2| - C_1 & \leq  G^+(F(z,w)) \\
& =  d^{1/n} G^+(z,w) \\
& \leq d^{1/n}( \max \{\log^+|z|, \log^+|w|\} + C_2).
\end{align*}
Exponentiating and using  $|F_1| \leq |F_2|+R$, we obtain
$$ |F(z,w)| \leq C \max\{ (|z|+1)^{d^{(1/n)}}, (|w|+1)^{d^{(1/n)}} \}. $$

Similarly, if $F(z,w) \in V^+$, then 
\begin{align*}
\log^+|F_1| - C_1 & \leq  G^-(F(z,w)) \\
& =  1/d^{(1/n)} G^-(z,w) \\
& \leq  1/d^{(1/n)} (\max \{\log^+|z|, \log^+|w|\} + C_2).
\end{align*}
Exponentiating and using  $|F_2| \leq |F_1|$, we obtain
$$ |F(z,w)| \leq C \max\{ (|z|+1)^{1/d^{(1/n)}}, (|w|+1)^{1/d^{(1/n)}} \}. $$

Hence $F$ has polynomial growth throughout $\CC^2$, hence must be a
polynomial.  

\bigskip
We show next that $F \in {\cal G}$.  Let $\tau(z,w) := (w,z)$.  By
\cite{fm}, we can write $H = \tau e_1 \cdots \tau e_m$, for some
elementary maps $e_j$.  Since $F^n=H$, $F$ must be a reduced word with
length at least 2.  There are four possibilities for the form of $F$.
The first is
$$ F = a_1 e_1' \cdots a_l e_l' $$
for some affine, non-elementary maps $a_j$ and some
elementary, non-affine maps $e_j'$.   By \cite{fm} or \cite{ar}, each
$a_j$ can be written $a_j = 
a_j^1 \tau a_j^2$, where $a_j^1$ and $a_j^2$ are affine and
elementary, and $a_1^1$ has the form
$$ a_1^1(z,w) = (bz +c w, w).  $$

Now, since $a_j^k$ is elementary, $F^n$ has the form $a_1^1 \circ
(p,q)$, where $p$ and $q$ are polynomials and $(p,q) = \tau e_1''
\cdots \tau e_{nl}''$.  By 
\cite{fm}, we have $\deg(q) > \deg(p)$, and likewise the degree of the
second coordinate function of $H=F^n$ is larger than the degree of the
first coodinate function.  This implies that $c=0$, so $a_1^1(z,w) =
(bz, w)$.  Replacing $a_2^1$ by $\sigma a_2^1$, where $\sigma(z,w) =
(z,bw)$, we obtain 
$$ F=\tau E_1 \cdots \tau E_l.$$
Hence $F \in {\cal G}$.  

The second case is
$$ F= e_1' a_2 \cdots a_l e_l'.$$
In this case, we can replace each $a_j$ by $a_j^1 \tau a_j^2$ as
before, and hence relabeling, we may assume $F = e_1' \tau \cdots \tau
e_l'$.  But then $H= F^n = e_1' \tau \cdots \tau e_{nl}'$, which
implies that the degree of the first coordinate of $H$ is larger than
the degree of the second coordinate, which is impossible.  Hence $F$
cannot have this form.

The third case is 
$$ F = e_1' a_2 \cdots e_{l-1}' a_l.  $$
As before, we may relabel to assume that $F = e_1' \tau \cdots
e_{l-1}' \tau a_l^2$.  But then $H = F^n = e_1' \tau \cdots e_k' \tau
a_l^2$, but also $H = \tau e_1 \cdots \tau e_m$.  Hence 
$$ I = (F^n)^{-1} H = (a_l^2)^{-1} (\tau (e_k')^{-1} \cdots
(e_1')^{-1}) (\tau e_1 \cdots \tau e_m).$$
But then $I$ has been written as a reduced word, which is impossible
from \cite{fm}.  Thus $F$ cannot have this form.

In the final case, we have
$$ F = a_1 e_1' \cdots e_{l-1}' a_l. $$
Again we may relabel and collect terms and assume $H=F^n = a_1^1 \tau
e_1' \cdots e_{k-1}' \tau a_l^2. $  Since $a_l^2$ is linear, we can
use an argument like that in the first case to relabel and replace
$a_1^1$ by the identity.  Since $a_l^2$ is elementary, we can write 
$a_l^2(z,w) = (az +bw+c, dw +e)$ with $a \neq 0$.  Applying $\tau e_1'
\cdots e_{k-1}' \tau$ to this, we see that the homogeneous polynomial
of highest degree in $F^n$ depends on $z$.  But a simple
inductive argument shows that the corresponding polynomial for $H$ is
independent of $z$.  Hence $F$ cannot have this form.  

Thus $F \in {\cal G}$.  $\BOX$

\bigskip

\begin{remark}
In general, a map $H$ can have distinct roots of a given order.  For
example, the map $H$ given by squaring $F(z,w)=(w, z+w^2)$ has
three square roots.  This is true because $(F \circ s)^2 =H$
for $s(z,w) = (\omega^2 z, \omega w)$, where $\omega^3=1$.  In fact,
one can check that these are the only possible square roots of $H$.
\end{remark}

\section{Roots of elementary maps}
In \cite{bf}, we showed that no H{\'e}non map can be the time-1 map of the
flow of a holomorphic vector field and gave a precise classification
of those maps which can be the time-1 map of such a flow.  In
\cite{forst} and \cite{afv}, it was shown that any flow of a
holomorphic vector field whose time-1 map is an elementary map is in
fact conjugate to a flow which is polynomial for all time.

In this section, we consider the set of elementary maps which are not
the time-1 map of any holomorphic flow and show that such maps have
roots of arbitrarily high order but that any root is conjugate to a
polynomial map.

The elementary maps which cannot be the time-1 map of a flow have the
form 
$$F(z,w) = (\beta^\mu(z+w^\mu q(w^r)), \beta w), $$
where $\beta$ is a primitive $r$th root of unity, $q(w) = w^k +
q_{k-1} w^k + \cdots q_1 w + q_0$, and $k \geq 1$.   A simple check
shows that if we replace $w^\mu q(w^r)$ by $(w^\mu q(w^r))/(lr+1)$ for
any $l \in \ZZ^+$, then the resulting map is an $(lr+1)$st root of
$F$.  

In general, maps of this form can have nonpolynomial roots.  For
instance, let $F(z,w) = (-(z+w(w^4+1)), -w)$ and let $k$ be any entire
function of one variable.  Define $\phi(z,w) = (i (z+w(w^4+1)/2 +w^3
k(w^4)), iw)$.  A simple check shows that $\phi^2 = F$, and $\phi$ is
nonpolynomial whenever $k$ is transcendental.

\bigskip
We claim that any root of $F$ is conjugate to a polynomial automorphism.
Suppose that $\phi$ is an automorphism of $\CC^2$ with $\phi^n =
F$.  Then $\phi F^r \phi^{-1} = F^r$, so an argument like that in
\cite[theorem 6.10]{fm} shows that $\phi(z,w) = (e^{g(w)} z +h(w), aw
+ b)$ for some 
entire $g$, $h$, and some $a, b \in \CC$, $a \neq 0$.  Since $\phi^n =
F$, we see that $a^n = \beta$ and $b(a^n-1)/(a-1) = 0$, so that $b=0$.

Using this form for $\phi$ and the fact that $\phi F^r = F^r \phi$, it
follows that $e^{g(w)}$ is a nonzero rational function, hence is a
constant, $c \neq 0$.  Moreover, since $\phi F = F \phi$, we see that
$c = a^\mu$.  Thus $\phi(z,w) = (a^\mu z + h(w), a w)$.  

Now, since $\phi^{rn} = F^r$, it follows that $\sum_{j=0}^{rn-1}
(a^\mu)^{-j} h(a^j w) = r w^\mu q(w^r)$.  Write $h = h_1 + h_2$, where
$h_2 = w^{\mu + kr} \tilde{h}_2$ for some entire $\tilde{h}_2$.  Then
the sum just given is valid with $h_2$ in place of $h$ and $0$ in place of
$r w^\mu q(w^r)$.  Hence by \cite{forst}, there exists $f$ entire such
that $f(a w) - a^\mu f(w) = h_2(w)$.  

A simple check shows that if $\psi(z,w) = (z + f(w), w)$, then
$\psi^{-1} \phi \psi$ is a polynomial, and in fact, $\psi^{-1} \phi^n
\psi = F$.  Thus any root of $F$ is conjugate to a polynomial
automorphism.  

\bigskip
Given any elementary map $F$ and a root $\phi^n=F$, one can ask if
$\phi$ is conjugate to a polynomial automorphism.  There are a few cases
such as the above where this result is relatively straightforward, but
in general, this seems to be a hard question.  For some results along
these lines in the case $F = I$, see \cite{ar}.

\bigskip

\small
\noindent Gregery T. Buzzard \\
Department of Mathematics\\
The University of Michigan\\
Ann Arbor, MI 48109, USA\\
and\\
MSRI \\
1000 Centennial Drive\\
Berkeley, CA 94720, USA\\

\noindent John Erik Forn\ae ss\\
Department of Mathematics\\
The University of Michigan\\
Ann Arbor, MI 48109, USA\\

\end{document}